\documentclass{article}
\usepackage{amsfonts}
\usepackage{amsmath}

\newtheorem{theorem}{Theorem}

\newtheorem{definition}[theorem]{Definition}

\input{tcilatex}

\begin{document}

\Large
\begin{center}
{\Large Mean Periodic Functions in Banach Spaces}

{\large \bigskip }

{\large K.\ Gowri Navada}
\end{center}

\textit{Abstract}{\normalsize : Here we give a necessary and
sufficient condition for a Banach space to be separable.}

{\large \bigskip }\textit{Key words:} mean periodic function,
vector valued measure

\bigskip

{\large For a Banach space }$B${\large \ over the complex field }$\mathbb{C}$%
{\large , let }$C(\mathbb{R},B)${\large \ denote the set of all
continuous
functions defined on the real line }$\mathbb{R}${\large \ taking values in }$%
B${\large \ with the compact convergence topology. When }$B=\mathbb{C}$%
{\large \ we write }$C(\mathbb{R})${\large \ for }$C(\mathbb{R},\mathbb{C})$%
{\large . For a function }$\phi ${\large \ in
}$C(\mathbb{R},B)${\large \
let }$\tau (\phi )${\large \ denote the closure in }$C(\mathbb{R},B)${\large %
\ of the span of all translates of }$\phi ${\large .}

\begin{definition}
{\large A function }$\phi ${\large \ in }$C(\mathbb{R},B)${\large
\ is said
to be\ mean periodic if }$\ ${\large \ }$\tau (\phi )\neq C(\mathbb{R},B)$%
{\large .}
\end{definition}

{\large We prove the following theorem.}

\begin{description}
\item[Theorem] {\large :}{\Large \ }{\large A Banach space }$B${\large \ is
separable if and only if not all functions in
}$C(\mathbb{R},B)${\large \ are mean periodic.}

\item[Preliminaries:] {\large Let }$B^{\ast \text{ }}${\large be the dual of
}$B${\large \ and let \ }$B_{s}^{\ast }${\large \ denote the space
}$B^{\ast \text{ }}${\large \ along with the weak}$^{\ast
}${\large \ topology on it. Let }$\mu ${\large \ be a countably
additive,}$\strut \mathstrut ${\large \ regular, vector valued,
Borel measure on }$\mathbb{R}${\large \ taking
values in }$B_{s}^{\ast }${\large . For a Borel set }$b^{\prime }\subseteq $%
$\mathbb{R}${\large , set }$\ \ \ \ \ \Vert \mu (b^{\prime })\Vert
=\sup_{\Vert v\Vert \leq 1}|<v,\mu (b^{\prime })>|${\large . This
is the norm in the Banach space }$B^{\ast \text{ }}.${\large \ For
\ a Borel set }$b ${\large \ \ in }$\mathbb{R}${\large \ let
}$P(b)${\large \ be
the set of all finite Borel partitions of }$b${\large . If }$%
\sup_{P(b)}\sum_{b^{^{^{\prime }}}}\Vert \mu (b^{\prime })\Vert <\infty $%
{\large \ (where the sum is taken over }$\{b^{\prime }\}$ {\large which form
a finite Borel partition of} $b${\large ),} {\large for all Borel sets }$b$%
{\large \ in }$B,${\large \ then }$\mu ${\large \ is said to be of bounded
variation.}
\end{description}

%\QTP{Body Math}
{\large Let }$M(\mathbb{R},B_{s}^{\ast })${\large \ be the space
of all vector valued Borel measures }$\mu ${\large \ on
}$\mathbb{R}${\large \ taking values in }$B_{s}^{\ast }${\large \
with the following properties:}

\begin{enumerate}
\item $\mu ${\large \ is countably additive and regular,}

\item $\mu ${\large \ has compact support in }$\mathbb{R}${\large ,}

\item $\mu ${\large \ is of bounded variation,}

\item {\large there exists a constant }$C>0${\large \ such that for any
Borel set }$b${\large \ of }$\mathbb{R}${\large , }$\Vert \mu
(b)\Vert \leq C ${\large \ where }$\Vert .\Vert ${\large \ denotes
the norm in the Banach space }$B^{\ast }${\large .}
\end{enumerate}

%\QTP{Body Math}
{\large Let }$\mu \in M(\mathbb{R},B_{s}^{\ast })${\large \ with support }$K$%
{\large . Let }$b_{1},b_{2},\cdots b_{n}${\large \ be a Borel partition of }$%
K${\large \ and }$v_{1},v_{2},\cdots v_{n}${\large \ be arbitrary
elements of }$B${\large . Then }$\phi =\sum_{i=1}^{n}v_{i}\chi
_{b_{i}}${\large \ is a simple function from }$\mathbb{R}${\large
\ to }$B${\large .\ Then }$\int \phi d\mu ${\large \ is defined
as}

\begin{center}
$\int \phi d\mu =\sum_{i=1}^{n}<v_{i},m(b_{i})>$
\end{center}

%\QTP{Body Math}
{\large Now let }$\phi :\mathbb{R}\longrightarrow B${\large \ be
any continuous function and let }$\mu ${\large \ and }$K${\large \
be as before.
For any }$\epsilon >0${\large \ we can get a finite Borel partition }$%
\{b_{i}\}${\large \ of }$K${\large \ such that }$\Vert \phi (x)-\phi
(y)\Vert <\epsilon ${\large \ for all }$x,y${\large \ in }$b_{i\text{ }}$%
{\large for all }$i${\large . Choose any point }$x_{i}\in
b_{i}${\large . Define the simple function }$\phi _{\epsilon
}=\sum \phi (x_{i})\chi _{b_{i}} ${\large . Then }$\int \phi d\mu
=\lim_{\epsilon \rightarrow 0}\int \phi _{\epsilon }d\mu ${\large
. This limit exists and is independent of the choice of the Borel
partition }$\{b_{i}\}${\large \ and also of the choice of the
points }$\{x_{i}\}${\large \ (see [S]).}

%\QTP{Body Math}
{\large For a function }$f\in C(\mathbb{R})${\large \ and a
measure }$\mu \in M(\mathbb{R},B_{s}^{\ast })${\large \ we define
the integral }$\int
fd\mu ${\large \ in a similar way: for a simple function }$%
f=\sum_{i}y_{i}\chi _{b_{i}}${\large \ where }$y_{i}\in
\mathbb{C}${\large \ and }$\{b_{i}\}${\large \ is a Borel
partition of }$K${\large , the support of }$\mu ${\large , define
}$\int fd\mu =\sum_{i}y_{i}\mu (b_{i})\in B_{s}^{\ast }${\large .
For }$f\in C(\mathbb{R})${\large , }$\int fd\mu
=\lim_{\epsilon \rightarrow 0}\int f_{\epsilon }d\mu ${\large , where }$%
f_{\epsilon }${\large \ are simple functions defined as before. Therefore }$%
\int fd\mu ${\large \ is an element of }$B_{s}^{\ast }${\large .}

%\QTP{Body Math}
{\large Singer's Theorem (see [S]): The dual of
}$C(\mathbb{R},B)${\large \ is the space
}$M(\mathbb{R},B_{s}^{\ast })${\large .}

%\QTP{Body Math}
{\large Using Hahn Banach Theorem it is easy to see that a function }$\phi $%
{\large \ is mean periodic if and only if there exists a nonzero measure }$%
\mu ${\large \ in }$M(\mathbb{R},B_{s}^{\ast })${\large \ such
that }$\phi \ast \mu =0${\large .}

\begin{description}
\item[Proof of the theorem:] {\large \ Suppose }$B${\large \ is not
separable. Let }$\phi ${\large \ be any function in }$C(\mathbb{R},B)$%
{\large . Since }$\mathbb{R}${\large \ is separable, }$\phi (\mathbb{R})$%
{\large \ is separable. Let }$B_{1\text{ }}${\large be the closure of the
space generated by }$\phi (\mathbb{R})${\large . It is easy to see that }$%
B_{1\text{ }}${\large is also separable. Moreover, range of any function in }%
$\tau (\phi )${\large \ is contained in }$B_{1\text{ }}.${\large \ Since }$%
B_{1\text{ }}${\large is a proper subspace of }$B${\large \ we can find a
function in }$C(\mathbb{R},B)${\large \ whose range is not contained in }$%
B_{1\text{. }}${\large Hence }$\tau (\phi )\neq
C(\mathbb{R},B)${\large .}
\end{description}

{\large Conversely, if }$B${\large \ is separable we will show
that there exists a function in }$C(\mathbb{R},B)${\large \ which
is not mean periodic. First we will take }$B${\large \ to be
}$l^{1}${\large , the space of
absolutely summable sequences in }$\mathbb{C}${\large . For each }$n\in N$%
{\large , let }$e_{n}=(0,\cdots 1,0,\cdots )${\large \ be the sequence in }$%
l^{1\text{ }}${\large where 1 occurs only at the }$n${\large th entry. For
each }$n\in N${\large \ we define an element }$\{a_{j}(n)\}_{j=1}^{\infty }$%
{\large \ in }$l^{1\text{ }}${\large and a sequence of real numbers }$%
\{\lambda _{j}(n)\}_{j=1}^{\infty }${\large \ such that:}

\begin{enumerate}
\item $a_{j}(n)${\large \ is nonzero for all }$j${\large \ and }$n${\large .}

\item {\large If we denote }$\sum_{j=1}^{\infty }|a_{j}(n)|=a(n),${\large \
then }$\sum_{n=1}^{\infty }a(n)<\infty ${\large \ (that is, }$%
\{a(n)\}_{n=1}^{\infty }${\large \ is also in }$l^{1}.$

\item $\lambda _{j}(n)${\large \ are all different. that is, if }$\lambda
_{j}(n)=\lambda _{i}(m)${\large \ then }$j=i${\large \ and }$n=m.$

\item {\large The sequence }$\{\lambda _{j}(n)\}_{j=1}^{\infty }${\large \
converges to a real number, say }$\lambda (n)${\large .}
\end{enumerate}

%\QTP{Body Math}
{\large For a real number }$s\in \mathbb{R}${\large \ define }$%
f_{n}(s)=\sum_{j=1}^{\infty }a_{j}(n)e^{i\lambda _{j}(n)s}.${\large \ Then }$%
f_{n}:\mathbb{R}\rightarrow \mathbb{C}${\large \ is continuous.
Put }$\phi (s)=${\large \ }$\sum_{n=1}^{\infty }${\large \
}$f_{n}(s)e_{n}${\large . It is easy to see that }$\phi
:\mathbb{R}\rightarrow l^{1}${\large \ is
continuous. The series }$\sum_{n=1}^{\infty }${\large \ }$f_{n}e_{n}${\large %
\ converges to the function }$\phi ${\large \ uniformly on }$\mathbb{R}$%
{\large . We will show that for this function }$\phi \in
C(\mathbb{R},l^{1}), ${\large \ }$\tau (\phi
)=C(\mathbb{R},l^{1})${\large .}

%\QTP{Body Math}
{\large Suppose not. Then by Singer's theorem \ there exists a
nonzero measure }$\mu \in M(\mathbb{R},l_{s}^{\infty })${\large \
such that }$\phi \ast \mu =0${\large . Let }$K$ {\large be the
compact support of }$\mu .$

%\QTP{Body Math}
\begin{equation*}
{\large 0=\phi \ast \mu =}\sum_{n=1}^{\infty }{\large (f}_{n}{\large e}_{n}%
{\large \star \mu )=}\sum_{n=1}^{\infty }{\large <e}_{n}{\large ,f}_{n}%
{\large \star \mu >}
\end{equation*}

%\QTP{Body Math}
{\large For the last equality refer to theorem III.2.6 in [S]. For }$n\in N$%
{\large \ and a Borel set }$b${\large \ in }$\mathbb{R}${\large ,
put }$\mu _{n}(b)=<e_{n},\mu (b)>${\large . Then }$\mu
_{n}${\large \ is a countably
additive, regular complex measure with compact support contained in }$K$%
{\large \ (see the proof of Singer's theorem in [S]). We can write }$\mu
(b)=(\mu _{1}(b),\mu _{2}(b),\cdots )\in l_{s}^{\infty }.${\large \ For any }%
$f\in C(\mathbb{R})${\large \ we have }$<e_{n},${\large \ }$\int fd\mu >=$%
{\large \ }$\int fd\mu _{n}${\large \ : this is true for
characteristic functions so also for simple functions. By the
limiting process one can prove it for any continuous function on
}$\mathbb{R}${\large . So by (1) we get }$\sum ${\large \
}$f_{n}\star \mu _{n}=0${\large . That is,}

%\QTP{Body Math}
\begin{equation*}
\sum_{n=1}^{\infty }\sum_{j=1}^{\infty }{\large a}_{j}{\large (n)}\overset{%
\symbol{94}}{\mu _{n}}{\large (\lambda }_{j}{\large
(n))e}^{i\lambda _{j}(n)s_{0}}{\large =0\forall s}_{0}{\large \in
\mathbb{R}}
\end{equation*}

%\QTP{Body Math}
{\large Since }$\overset{\symbol{94}}{\mu _{n}}(\lambda _{j}(n))${\large \
are uniformly bounded by }$\Vert \mu \Vert ${\large \ the left hand side is
an almost periodic function. Because }$a_{j}(n)\neq 0${\large \ for all }$%
j,n ${\large , this implies }$\overset{\symbol{94}}{\mu _{n}}(\lambda
_{j}(n))=0. ${\large \ The zero's of the holomorphic function }$\overset{%
\symbol{94}}{\mu }_{n}${\large have a limit point }$\lambda (n).${\large \
Therefore }$\mu _{n}=0${\large \ for all }$n${\large \ and hence }$\mu =0$%
{\large , a contradiction.}

%\QTP{Body Math}
{\large Now let }$B${\large \ be an arbitrary separable Banach space. Since }%
$B${\large \ is separable we can find a countable set }$\{h_{n}\}${\large \
in }$B${\large \ such that }$\Vert h_{n}\Vert =1${\large \ for all }$n$%
{\large \ and the subspace }$H${\large \ generated by }$\{h_{n}\}${\large \
is dense in }$B${\large . Define a function }$\psi :\mathbb{R}\rightarrow B$%
{\large \ by }$\psi (s)=\sum_{n=1}^{\infty }f_{n}(s)h_{n}${\large \ where }$%
f_{n}${\large 's are defined as before. Then this series converges in }$B$%
{\large \ and the function }$\psi ${\large \ is continuous. Let }$g\in $%
{\large \ }$C(\mathbb{R})${\large \ be any function. Then we claim that }$%
gh_{1}\in \tau (\psi )${\large :}

%\QTP{Body Math}
{\large Since }$\tau (\phi )=C(\mathbb{R},l^{1})${\large \ we know that }$%
ge_{1\text{ }}${\large is a limit (in
}$C(\mathbb{R},l^{1})${\large ) of a
sequence of linear combinations of translates of }$\phi ${\large . That is, }%
$ge_{1}=\lim_{m\rightarrow \infty }\Phi _{m}${\large \ where for each }$m\in
N${\large , }$\ \ \Phi _{m}${\large \ is a finite sum}$\ \Phi _{m}=\sum
c_{i}\phi _{y_{i}}${\large . Define }$\Psi _{m}=\sum c_{i}\psi _{y_{i}}$%
{\large \ for each }$m.${\large Then }$gh_{1}=\lim_{m\rightarrow
\infty }\Psi _{m}${\large : Let }$\epsilon >0${\large \ and
}$C${\large \ be a compact subset of }$\mathbb{R}${\large . For
any }$s\in C$

%\QTP{Body Math}
\begin{eqnarray*}
\Vert g(s)h_{1}-\Psi _{m}(s)\Vert &=&\Vert g(s)h_{1}-\sum c_{i}\psi
_{y_{i}}(s)\Vert \\
&=&\Vert g(s)h_{1}-\sum_{n=1}^{\infty
}\sum_{i}c_{i}(f_{n})_{y_{i}}(s)h_{n}\Vert \\
&\leq &|g(s)-\sum_{i}c_{i}(f_{1})_{y_{i}}(s)+\sum_{n\geq
2}|\sum_{i}c_{i}(f_{n})_{y_{i}}(s)| \\
&=&\Vert g(s)e_{1-}\Phi _{m}(s)\Vert _{l^{1}}\leq \epsilon
\end{eqnarray*}

%\QTP{Body Math}
{\large for all }$m${\large \ sufficiently large.}

%\QTP{Body Math}
{\large Similarly we can prove that any finite sum }$\sum g_{i}h_{i}${\large %
\ in }$C(\mathbb{R})\otimes H${\large \ \ is in }$\tau (\psi
).${\large \
Since }$\tau (\psi )${\large \ is closed and }$C(\mathbb{R})\otimes H$%
{\large \ is dense in }$C(\mathbb{R})\otimes B${\large , we get }$C(\mathbb{R%
})\otimes B\subseteq ${\large \ }$\tau (\psi ).${\large \ But }$C(\mathbb{R}%
)\otimes B${\large \ is dense in }$C(\mathbb{R},B)${\large \ (see
[S]). Hence }$\tau (\psi )=C(\mathbb{R},B)${\large .}

%\QTP{Body Math}
$\medskip $

%\QTP{Body Math}
{\LARGE References}

%\QTP{Body Math}
${\LARGE \medskip }$

%\QTP{Body Math}
{\large [K] J.P. Kahane, Lectures on Mean Periodic Functions; Tata Institute
of Fundamental Research, Mumbai, 1959.}

%\QTP{Body Math}
{\large [S] J.Schmets, Spaces of Vector Valued Continuous Functions; Lecture
Notes in Mathematics 1003, Springer Verlag, Berlin, 1983.}

%1\QTP{Body Math}
{\large [T] F.Treves, Topological Vector Spaces,
Distributions and Kernels; Academic Press, New York, 1967.}

\end{document}